\documentclass{amsart}
\usepackage{amsmath,amscd,amsopn,amsthm,amsxtra,upref,amssymb,latexsym}
\theoremstyle{plain}
\newtheorem{theorem}{Theorem}

\theoremstyle{definition}

\theoremstyle{remark}

\numberwithin{equation}{section}
\newcommand{\al}{\alpha}

\newcommand{\R}{\mathbb R}
\newcommand{\N}{\mathbb N}

\newcommand{\C}{\mathbb C}
\newcommand{\Sn}{\mathbb S^{n-1}}
\newcommand{\Rn}{\mathbb R^n}

\newcommand{\Rm}{\mathbb R^{n+1}}
\begin{document}

\title{The Taylor series of the Gaussian kernel}
%%%%%%%%%%%%%%%%%%
%Author information
\author{L. Escauriaza}
\address[L. Escauriaza]{Universidad del Pa{\'\i}s Vasco / Euskal Herriko
Unibertsitatea\\ Departamento de Matematicas\\ Apartado 644, 48080 Bilbao, Spain.} 
\email{luis.escauriaza@ehu.es}
\thanks{\thanks{The author is supported  by MEC grant MTM2004-03029 and by the UPV/EHU grant 9/UPV 00127.310-15969/2004}}
%%%%%%%%%%%%%%%%%%%
%\date{}
%%%%%%%%%%%%%%
\begin{abstract}
We describe a formula for the Taylor series expansion of the Gaussian kernel around the origin of $\Rn\times\R$.
\end{abstract}
\dedicatory{"From some people one can learn more than mathematics"}
\maketitle

\section{Introduction}\label{S:1}
The explicit formulae for the power series expansion at the origin of the fundamental solution of the Laplace operator in $\Rn$, $n\ge 2$, are well known. In particular, when
\begin{equation*}
\Gamma(x,y)=
\frac 1{\omega_n(2-n)}|x-y|^{2-n}\ ,
\end{equation*}
$\omega_n$ is the surface measure of the unit sphere $\Sn$ in $\R^n$ and $n>2$, or when
\begin{equation*}
\Gamma(x,y)=\frac 1{2\pi}\log{|x-y|}
\end{equation*}
and  $n=2$, the following hold:
\begin{equation}\label{E:Laplace1}
\varphi(x)=\int_{\Rn}\Gamma(x,y)\triangle\varphi (y)\, dy\quad ,\text{when}\  \varphi\in C_{0}^{\infty}(\Rn)\ ,
\end{equation}
\begin{equation}\label{E:Laplace2}
\Gamma(x,y)=
\sum_{k=0}^{+\infty}\frac{|x|^k}{|y|^{k+n-2}}Z_k(x'\cdot y')\ ,\ \text{when}\ 0\le |x|<|y|\ ,
\end{equation}
\begin{equation}\label{E:Laplace3}
Z_k(x'\cdot y')=\begin{cases}
-\frac 1{2k+n-2}Z_{x'}^{(k)}(y')\ &,\ \text{when}\  2k+n-2>0\ ,
\\
\frac 1{2\pi}\log{|y|}\ &,\ \text{when}\  2k+n-2=0\ .
\end{cases}
\end{equation}
Here, $x'=x/|x|$ and $Z_{x'}^{(k)}(y«)$ is the \emph{zonal harmonic of degree k}, i.e., the kernel  of the projection operator of $L^2(S^{n-1})$ onto the spherical harmonics of degree $k\ge 0$. We recall that the  the spherical harmonics of degree $k\ge 0$ are the eigenfunctions of the spherical Laplacian on $\Sn$ and corresponding to the eigenvalue, $k(k+n-2)$. See  \cite[(2.17)]{GilbargTrudinger} for \eqref{E:Laplace1} and \cite[Chapter IV]{Stein1} for \eqref{E:Laplace2} and \eqref{E:Laplace3}.
\footnote{A simple way to derive \eqref{E:Laplace2} and \eqref{E:Laplace3} is to solve, $\triangle u=f$ in $\Rn$, via the method of separation of variables in spherical coordinates and then, to compare the solution, which the latter method yields, with the one obtained via the convolution with the fundamental solution.}

For fixed $z'$ in $S^{n-1}$ and $k\ge 0$, the function
\[E_{k}(x)=|x|^kZ_k(x'\cdot z')\]
is a homogeneous harmonic polynomial of degree $k$, i.e., 
\[\triangle E_{k}=0\ \text{and}\  E_{k}(\lambda x)=\lambda^{k}E_{k}(x)\ ,\]
 when $\lambda\ge 0$ and $x$ is in $\Rn$. The function, \[|y|^{2-n-k}Z_k(z'\cdot y')\] is harmonic in $\Rn\setminus\{0\}$ and is  homogeneous of degree, $2-k-n$. In fact, it is is the Kelvin transformation of $E_k$. 
Here recall, that the Kelvin transformation $v$ of a  function $u$ is 
\begin{equation*}
v(x)=|x|^{2-n}u(x/|x|^2)\quad\text{,}\quad \triangle v(x)=|x|^{-n-2}\triangle u(x/|x|^2)\ ,
\end{equation*}
\begin{equation*}
\varphi=\sum_{k=0}^{+\infty}Z_k(\varphi)\quad \text{,}\quad \|\varphi\|_{L^2(\Sn)}^2=\sum_{k=0}^{+\infty}\|Z_k(\varphi)\|_{L^2(\Sn)}^2\ ,
\end{equation*}
when $\varphi$ is in $C^\infty(\Sn)$ and where
\[Z_k(\varphi)=\int_{\Sn}Z_{x'}^{(k)}(y')\varphi(y')\,dy'\ \text{and}\ k\ge 0\ .\]
Moreover, the formula in \eqref{E:Laplace2} is  the Taylor series expansion of $\Gamma(x,y)$ around $x=0$, for each fixed $y\neq 0$ in $\R^n$. (See \cite[Chapter IV]{Stein1}).

The Taylor series \eqref{E:Laplace2} contains relevant information, which has had important applications in the mathematics of the last century and among others it has shown to be useful to obtain estimates leading to sharp results of strong and weak unique continuation for elliptic operators on $\Rn$. This can be seen in \cite{jk}, \cite{w1} and \cite{kt}.

The Gaussian kernel
\begin{equation*}\label{E:forwardGaussian}
G(x,t,y,s)=\begin{cases}
\left(4\pi (t-s)\right)^{-n/2}e^{-|x-y|^2/4(t-s)}\ &,\ \text{when} \ s<t\ ,
\\0\ &,\ \text{when} \ s>t\ ,\end{cases}
\end{equation*}
is the fundamental solution of the heat operator in $\Rm$, i.e.
\begin{equation*}
f(x,t)=-\int_{-\infty}^{t}\int_{\Rn}G(x,t,y,s)(\triangle f-\partial_{s}f)\, dyds\ ,\ \text{when}\ f\in C_{0}^{\infty}\left(\Rm\right)\ .
\end{equation*}

As far as the author knows (and this is  rather  surprising), it seems that nobody has written down and publish an explicit formula for the Taylor series expansion of $G(x,t,y,s)$ around the origin of $\Rm$, when $(y,s)$ in $\Rm$, $s<0$, is fixed. The purpose of this note is to fill in this gap.

To simplify the notation, we choose to give the formula for the Taylor series expansion of the fundamental solution of the backward heat equation,
\begin{equation}\label{E:backwardGaussian}
G_{b}(x,t,y,s)=\begin{cases}
\left(4\pi (s-t)\right)^{-n/2}e^{-|x-y|^2/4(s-t)}\ &,\ \text{when} \ t<s\ ,
\\0\ &,\ \text{when} \ t>s\ ,\end{cases}
\end{equation}
when $s$ is positive and $(y,s)$ in $\Rm$ is fixed. The Taylor series for the Gaussian kernel follows from the identity
\begin{equation*}
G(x,t,y,s)=G_{b}(x,-t,y,-s)\ .
\end{equation*}

The Hermite functions, $h_k$, are defined as
\begin{equation*}
h_k(x)=\left(2^kk!\sqrt\pi\right)^{-\frac 12}(-1)^ke^{x^2/2}\frac{d^k}{dx^k}\left(e^{-x^2}\right)\ ,\ k\ge 0\ , \ x\in\R
\end{equation*}
and $h_k=H_k(x)e^{-x^2/2}$, where $H_k$ is a Hermite polynomial of degree $k$. 

The Hermite functions on $\Rn$, $\phi_\al$, $\al=(\al_1,\dots,\al_n)$ in ${\N}^n$, are the product of 
the one-dimensional Hermite functions $h_{\al_j}$, $j=1,\dots,n$
\begin{equation*}
\phi_\al (x)=\prod_{j=1}^n h_{\al_j}(x_j)\ .
\end{equation*}
They form a complete orthonormal system in $L^2(\Rn)$ and if, $H=\triangle -|x|^2$, is the Hermite
operator, $H\phi_{\al}=-\left(2\left|\al\right|+n\right)\phi_{\al}$, where $|\al|=\al_1+\dots+\al_n$. The kernel
\begin{equation*}\label{E:hermite kernels}
\Phi_k(x,y)=\sum_{|\al|=k}\phi_\al(x)\phi_{\al}(y)
\end{equation*}
is the kernel of the projection operator of $L^2(\Rn)$ onto the Hermite functions of degree, $k\ge 0$, and when  $\varphi\in C_0^\infty(\Rn)$,
\begin{equation*}
\varphi=\sum_{k=0}^{+\infty}P_k(\varphi)\quad\text{,}\quad \|\varphi\|_{L^2(\Rn)}^2=\sum_{k=0}^{+\infty}\|P_k(\varphi)\|_{L^2(\Rn)}^2
\end{equation*}
and where
\begin{equation*}
P_k(\varphi)=\int_{\Rn}\Phi_k(x,y)\varphi(y)\,dy\ ,\ k\ge 0\ .
\end{equation*}

The reader can find the proofs of the latter results in \cite[Chapter 1]{Thangavelu}.

What seems to be the counterpart of the Kelvin transformation in the parabolic setting is the Appell transformation $v$ of a function $u$ (\cite{Appell}, \cite[pp. 282]{Doob}):
\begin{equation*}
v(x,t)=|t|^{-n/2}e^{-|x|^2/4t}u(x/t ,1/t)
\end{equation*}
and
\begin{equation*}
\triangle v-\partial_tv=|t|^{-2-n/2}e^{-|x|^2/4t}(\triangle u+\partial_tu)(x/t ,1/t)\ .
\end{equation*}

The Appell transformation maps backward caloric functions into forward caloric functions.

A calculation shows that 
\begin{equation*}
Q_\al(x,t)=t^{k/2}\phi_\al(x/2\sqrt t)e^{|x|^2/8t}\ , \ \al\in\N^n\ ,\  |\al| =k\ ,
\end{equation*}
 is backward caloric and in fact, it is a backward caloric polynomial in the $(x,t)$-variables, which is homogeneous of degree $k=|\al|$ in the parabolic sense, i.e.,
\begin{equation*}
\triangle Q_\al+\partial_t Q_\al=0\ \text{in}\  \Rm\ \text{and}\  Q_\al(\lambda x,\lambda^2t)=\lambda^kQ_\al(x,t)\ ,
\end{equation*}
when $\lambda\ge 0$ and $(x,t)$ is in $\Rm$(The later follows because a Hermite polynomial $H_k$ is an even function, when $k$ is even and an odd function, when $k$ is odd). At the same time and in analogy with the what happens with the Kelvin transformation of the harmonic function, $|x|^{k}Z_k(z'\cdot x')$, the function
\[s^{-(k+n)/2}\phi_\al(y/2\sqrt s)e^{-|y|^2/8s}\ ,\] 
is forward caloric and is the Appell transformation of $Q_\al$.

Having gathered all this data, it is possible to describe and write down the Taylor series expansion of the backward Gaussian kernel \eqref{E:backwardGaussian} at the origin of $\Rm$. We do it in the following theorem:
\begin{theorem}\label{T:1} The following identity holds, when $ t<s$, $s>0$, and $x$, $y$ are in $\Rn$
\begin{equation}
G_{b}(x,t,y,s)=(4s)^{-n/2}e^{|x|^2/8t}\left(\sum_{k=0}^{+\infty}\left(t/s\right)^{k/2}\Phi_k(x/2\sqrt t,y/2\sqrt s)\right)e^{-|y|^2/8s}\ .
\end{equation}
\end{theorem} 

The proof of Theorem \ref{T:1} is given in section \ref{S:2} and it follows from a well known identity: the \emph{generating formula} for the kernels, $\Phi_k$ (See \eqref{E:generating formula} below). 

The commentaries in \cite[pp. 582--583]{Stein2}, which are made with the purpose to explain the reader a simple approach to prove the identity \eqref{E:generating formula} and in particular, the reference \cite[pp. 335--336]{Feller}, show that \emph{Theorem \ref{T:1} was probably already known to some authors, though not explicitly written down and published.} In fact, the approach suggested in \cite[pp. 582--583]{Stein2} to prove the identity \eqref{E:generating formula}, follows precisely the inverse path of the one we follow in section 2 to prove Theorem \ref{T:1}. Thus, Theorem \ref{T:1} was  probably known by W. Feller and E.M. Stein.

In the same way as the formula for the Taylor series of the fundamental solution of the Laplace operator has been useful to derive results of unique continuation  for elliptic operators, the formula in Theorem \ref{T:1} is what, in a certain sense, is behind  the positive results of unique continuation for parabolic equations in \cite{e}, \cite{ev}, \cite{f} and \cite{ef}.

The argument is section \ref{S:2} gives a clue of how to proceed to find the Taylor series expansion of the fundamental solution of the Schr\"odinger operator, $\triangle +i\partial_{t}$, around the origin of $\Rm$ and the corresponding building pieces of the solutions of the Schr\"odinger equation: ``the \emph{Schr\"odinger} homogeneous polynomials of degree $k\ge 0$''. 

\section{Proof of Theorem \ref{T:1}}\label{S:2}
\begin{proof}
Recall the \emph{generating formula} for the projection kernels, $\Phi_{k}$ \cite{Thangavelu}
\begin{equation}\label{E:generating formula}
\sum_{k=0}^{+\infty}\Phi_k(x,y)\xi^k=\pi^{-\frac n2}\left(1-\xi^2\right)^{-\frac n2}e^{-\frac 12\frac{1+\xi^2}{1-\xi^2}\left(|x|^2+|y|^2\right)+
\frac{2\xi xy}{1-\xi^2}}\ ,\ \text{when}\quad |\xi |<1\ ,\ \xi\in\C\ .
\end{equation}
Replace $x$ by $x/2\sqrt t$, $y$ by $y/2\sqrt s$ and take $\xi=\sqrt{t/s}$ in \eqref{E:generating formula}, when $0\le t<s$. It gives
\begin{equation}\label{E:generating formula1}
s^{-n/2}\sum_{k=0}^{+\infty}\left(t/s\right)^{k/2}\Phi_k(x/2\sqrt t,y/2\sqrt s)=\pi^{-\frac n2}\left(s-t\right)^{-\frac n2}e^{-\frac{s+t}{s-t}\left(|x|^2/8t+|y|^2/8s\right)+
\frac{xy}{2(s-t)}}\ .
\end{equation}
Then, multiply \eqref{E:generating formula1}  by $4^{-n/2}e^{|x|^2/8t-|y|^2/8s}$ to get that the identity
\begin{multline*}\label{E:generating formula4}
(4s)^{-n/2}e^{|x|^2/8t}\left(\sum_{k=0}^{+\infty}\left(t/s\right)^{k/2}\Phi_k(x/2\sqrt t,y/2\sqrt s)\right)e^{-|y|^2/8s}\\=\left(4\pi(s-t)\right)^{-\frac n2}e^{-|x-y|^2/4(s-t)}
\end{multline*}
holds, when $0\le t<s$, $s>0$ and  $x$, $y$ are in $\Rn$, and Theorem \ref{T:1} follows.
\end{proof}

%%%%%%%%%%%%%%%%%%%%%%

%%%%%%%%%%%%%%%%%%%%%%%%%%

\begin{thebibliography}{99}

\bibitem{Appell} P. Appell, \emph{Sur l'\'equation $\partial^2z/\partial^2x-\partial z/\partial y=0$ et la th\'eorie du chaleur,} J. Math. Pures. Appl. (4) \textbf{8} (1892), 187--216.

\bibitem{Doob}
J.L. Doob, \emph{Classical potential theory and its probabilistic counterpart,} Springer-Verlag New-York Berlin, Heildeberg Tokyo, 1984.

\bibitem{e} L. Escauriaza, \emph{Carleman inequalities and the heat operator,} Duke Math.
J. \textbf{104}, n.1 (2000), 113-127.

\bibitem{ef} L. Escauriaza, F.J. Fern\'andez, \emph{Unique continuation for parabolic
operators,} Ark. Mat. \textbf{41} (2003), 35--60.

\bibitem{ev} L. Escauriaza, L. Vega, \emph{Carleman inequalities and the heat operator II,}
Indiana U. Math. J. \textbf{50}, n.3 (2001), 1149--1169.

\bibitem{Feller}
W. Feller, \emph{An introduction to probability theory and its applications, Vol. II,} Wiley.

\bibitem{f} F.J. Fern\'andez, \emph{Unique continuation for parabolic operators II,}
Comm. Part. Diff. Equat. \textbf{28} n. 9 \& 10 (2003), 1597--1604.

\bibitem{GilbargTrudinger}
D. Gilbarg, N.S. Trudinger, \emph{Elliptic partial differential equations of second order,} Springer-Verlag Berlin, Heildeberg, New-York, Tokyo, 1983.

\bibitem{jk} D. Jerison, C.E. Kenig, \emph{Unique continuation and absence of positive eigenvalues for Schr\"odinger operators,} Annals of Math. \textbf{121} (1985), 463--488.

\bibitem{kt} H. Koch, D. Tataru, \emph{Carleman estimates and unique continuation for second order elliptic equations with nonsmooth coefficients,} Comm. Pure Appl. Math.\textbf{54} n. 3 (2001), 339--360.

\bibitem{Stein1}
E.M. Stein, \emph{Introduction to Fourier analysis on euclidean spaces,} Princeton, New Jersey. Princeton University Press, 1975.

\bibitem{Stein2}
\bysame, \emph{Harmonic Analysis: Real-variable Methods, Orthogonality and Oscillatory Integrals,} Princeton University Press, 1993.  

\bibitem{Thangavelu}
S. Thangavelu, \emph{Lectures on Hermite and Laguerre expansions,} Princeton Univ. Press, Princeton, New Jersey, 1993.

\bibitem{w1} T. Wolff, \emph{Unique continuation for $|\triangle u|\le V|\nabla u|$ and related problems,} Rev. Mat. Iberoamericana \textbf{6} (1990), 155--200.

\end{thebibliography}
\end{document}